%% file: main.tex
\documentclass{article}
\usepackage[left=1in,right=1in,top=1in,bottom=1.3in]{geometry}
\usepackage{amsmath, amssymb, amsthm, latexsym, amsfonts, indentfirst, xcolor, mathtools, manfnt, microtype, enumitem, graphicx}
\usepackage[utf8]{inputenc}
\usepackage[english]{babel}
\usepackage{csquotes}
\graphicspath{ {./} }

\usepackage[breaklinks]{hyperref}
\usepackage{cleveref}
\hypersetup{
	colorlinks = true, 
	urlcolor = cyan, 
	linkcolor = teal, 
	citecolor = cyan 
}

\usepackage[backend=biber, sorting=nty, giveninits=true, maxbibnames=99]{biblatex}
\renewbibmacro{in:}{}
\DeclareFieldFormat{pages}{#1}

\usepackage{subfiles}

\newtheorem{Conjecture}{Conjecture} 
\newtheorem{Corollary}{Corollary} 
 
\newtheorem{Lemma}{Lemma} 
\newtheorem{Proposition}{Proposition} 
\newtheorem{Theorem}{Theorem}

\newcommand{\z}{{\mathbf z}}
\newcommand{\R}{{\mathbb R}}
\newcommand{\Z}{{\mathbb Z}}

\newcommand{\Q}{{\mathbb Q}}

\newcommand{\V}{{\mathbf v}}

\newcommand{\refheart}{$(\hyperlink{heart}{\heartsuit})$}
\newcommand{\refdiamond}{$(\hyperlink{diamond}{\diamondsuit})$}
\newcommand{\refclub}{$(\hyperlink{club}{\clubsuit})$}
\newcommand{\refspade}{$(\hyperlink{spade}{\spadesuit})$}
\newcommand{\refstar}{$(\hyperlink{star}{\bigstar})$}

\title{Monochromatic triangles in the max-norm plane}
\author{
    Alexander Natalchenko\thanks{MIPT, Moscow, Russia.  Email:~\href{mailto:natalchenko.ae@gmail.com}{\tt natalchenko.ae@gmail.com}}
    \and
    Arsenii Sagdeev\thanks{Alfréd Rényi Institute of Mathematics, Budapest, Hungary.  Email:~\href{mailto:sagdeevarsenii@gmail.com}{\tt sagdeevarsenii@gmail.com}.}
}

\date{}

\begin{document}

\maketitle

\begin{abstract}
For all non-degenerate triangles $T$, we determine the minimum number of colors needed to color the plane such that no max-norm isometric copy of $T$ is monochromatic.
\end{abstract}

\section{Introduction}

\subfile{introduction}

\section{Preliminaries}

\subfile{prelim}

\section{Proof of Theorem~\ref{main_theorem}}

\subfile{proof}

\section{Proof of Theorem~\ref{irrational_theorem}} \label{sect_irr}

\subfile{proof2}

\section{Concluding remarks}

\subfile{conclusions}

\printbibliography
\end{document}

%% file: introduction.tex
Modern \textit{combinatorial geometry} is both deep and wide with powerful tools galore. Nevertheless, some of its questions, that may look quite simple at first glance, repel all the attempts to answer them for several decades. Perhaps the most taunting problem of this sort is due to Nelson\footnote{See Soifer's investigation on the tangled origin of this question in~\cite[Section~3]{soifer2009}.} who asked in 1950 to find the \textit{chromatic number} $\chi(\R^2)$ of the Euclidean plane defined as the minimum number of colors needed to color the plane $\R^2$ such that no two points at unit Euclidean distance apart are of the same color. Despite the long history of research, it is only known that $5 \leq \chi(\mathbb{R}^2) \leq 7$, where the lower bound was obtained only in 2018, see~\cite{deGrey, exoo2020chromatic}. For multidimensional versions of this problem, see~\cite{CHERKASHIN2018125} and the references therein.

In their celebrated trilogy \cite{euclidian_ramsey_theorems_1, euclidian_ramsey_theorems_2, euclidian_ramsey_theorems_3}, Erdős, Graham, Montgomery, Rothschild, Spencer, and Straus laid the foundation of \textit{Euclidean Ramsey theory} which deals with questions of similar flavor but with more complex configurations forbidden to be monochromatic, see the survey~\cite{GrahamSurvey} by Graham. After a pair of points, the second simplest configuration is of course a triangle\footnote{Here and in what follows, we identify a triangle with the set of its three vertices.}. We denote by $\chi(\R^2, T)$ the minimum number of colors needed to color the plane such that no \textit{isometric} (i.e., translated and rotated) copy of a triangle $T$ is monochromatic. Erdős et al. conjectured in~\cite[Conjecture~3]{euclidian_ramsey_theorems_3} that  $\chi(\R^2, T) \ge 3$ for all triangles $T$ except for an equilateral one\footnote{For an equilateral triangle $\triangle$, the same group of authors observed that $\chi(\R^2, \triangle) = 2$ and conjectured in~\cite[Conjecture~1]{euclidian_ramsey_theorems_3} that the corresponding two-coloring is unique, which was later disproved in~\cite{jelinek2009monochromatic}.}, i.e., that two colors are never enough. Despite the efforts of various researchers, this conjecture was verified only for a few special families of triangles, see \cite{currier2024twocoloring, euclidian_ramsey_theorems_3, SHADER1976385, shkredov2015some}. From the other direction, it is easy to see that $\chi(\R^2, T) \le \chi(\R^2)$ and thus $\chi(\R^2, T) \le 7$ for all triangles $T$. Perhaps surprisingly, no better general upper bound is known, though Graham conjectured, see \cite[Conjecture 11.1.3]{GrahamSurvey} and~\cite{soifer1991triangles}, that $\chi(\R^2, T) \le 3$ for all triangles $T$, which was confirmed for `not very flat' triangles in~\cite{aichholzer2019triangles}. Let us also mention that currently the best bounds for multidimensional variant of this problem were recently obtained in~\cite{kupavskii2022cutting}. See also \cite{cheng2023euclidean, frankl2023monochromatic, fuhrer2024progressions, voronov2024chromatic} for related  geometric Ramsey results in the plane. 

\begin{figure}[h]
	\includegraphics[scale=0.19]{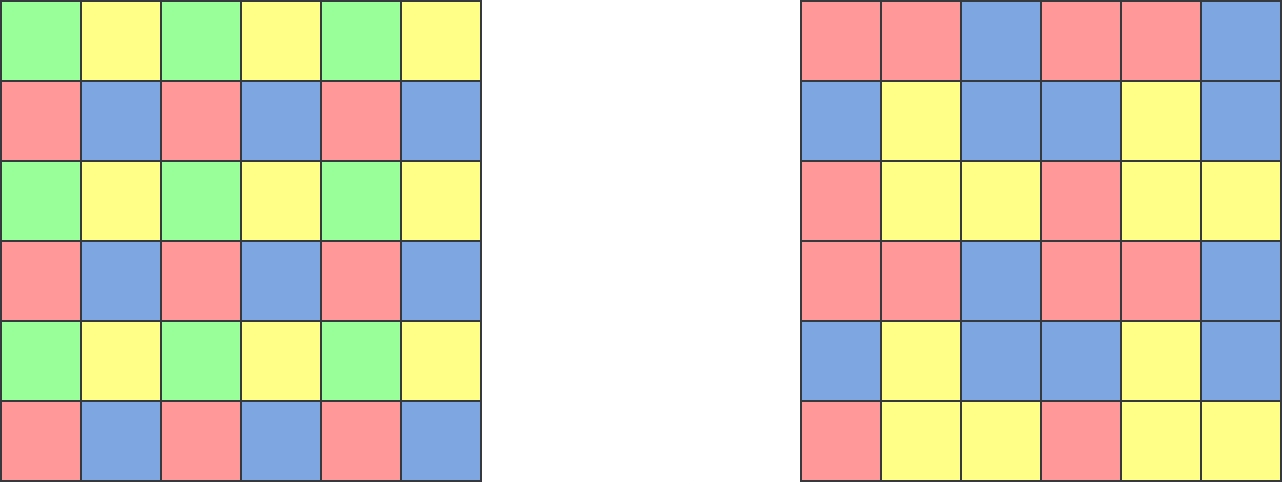}
	\centering
	\caption{A four-coloring of the plane with no monochromatic points at unit $\ell_\infty$-distance apart, and a three-coloring of the plane with no monochromatic $\ell_\infty$-isometric copies of degenerate triangles $T(a,b,a+b)$ for all coprime integers $a,b$ such that $a \equiv b \pmod{3}$. The side length of all the squares equals $1$.}
	\label{fig1}
\end{figure}

In this paper, we continue the line of research from \cite{frankl2021max, kirova2022two, kupavskii2021all, kupavskii2022spectrum} and consider a max-norm counterpart of the aforementioned problem. To give a formal definition, let us recall some basic notions and facts first. The \textit{$\ell_\infty$-distance} between $\z_1=(x_1,y_1)$, $\z_2=(x_2,y_2)\in \R^2$ is given by $\|\z_1-\z_2\|_\infty = \max\{|x_1-x_2|, |y_1-y_2|\}$. In contrast to the Euclidean case, it is easy to find the exact value of $\chi(\R_\infty^2)$ defined as the minimum number of colors needed to color the plane such that no two points at unit $\ell_\infty$-distance apart are of the same color: see the folklore proof that $\chi(\R_\infty^2)=4$ in the left-hand side of \Cref{fig1} or e.g. in~\cite[Section~2.1]{kupavskii2021all} for more details. A subset $T' \subset \R^2$ is called an \textit{$\ell_\infty$-isometric copy of $T \subset \R^2$}, if there exists a bijection $f:T \to T'$ such that $\|\z_1-\z_2\|_\infty = \|f(\z_1)-f(\z_2)\|_\infty$ for all $\z_1,\z_2 \in T$. Finally, we denote by $\chi(\R_\infty^2, T)$ the minimum number of colors needed to color the plane such that no $\ell_\infty$-isometric copy of $T$ is monochromatic. As earlier, it is easy to see that $\chi(\R_\infty^2, T) \le \chi(\R_\infty^2)$ and thus $\chi(\R_\infty^2, T)$ is equal to either $2$, or $3$, or $4$ for every triangle~$T$. We will show that all three of these options indeed take place.

Note that the value of $\chi(\R_\infty^2,T)$ depends only on the side lengths of $T$ and is independent of the particular position of $T$ in the plane. For all positive $a \leq b \leq c$ satisfying the triangle inequality $c \le a+b$, let $T(a, b, c)$ be an arbitrary triple of points in the plane with pairwise $\ell_\infty$-distances between them of $a, b, c$. This triangle is \textit{degenerate}\footnote{Note that three vertices of a degenerate triangle are not necessarily collinear, as the example of $(0,0), (1,0),$ and $(2,1)$ shows.} if $c=a+b$, otherwise it is \textit{non-degenerate}. Observe that if both fractions $\frac{a}{c}$ and $\frac{b}{c}$ are rational, then after a proper scaling, we can assume without loss of generality that $a,b,c$ are coprime integers. For all such non-degenerate triangles $T$, our next result gives the exact value of $\chi(\R_\infty^2,T)$.

\begin{Theorem} \label{main_theorem}
	Let $a \le b \le c$ be positive integers such that $c<a+b$ and $\gcd(a,b,c)=1$. Put $T=T(a,b,c)$. If \hypertarget{heart}{}
	\begin{itemize}[noitemsep,topsep=0pt]
		\item[$(\heartsuit)$] $a+b+c$ is odd, or \hypertarget{diamond}{}
		\item[$(\diamondsuit)$]  $a$ and $b$ are odd, $c\ge a+b-\gcd(a,b)$,
	\end{itemize}
	\noindent then $\chi(\R_\infty^2, T)=2$. Otherwise, $\chi(\mathbb{R}_\infty^2, T) = 3$.
\end{Theorem}

Our next complementary result covers the remaining case of `irrational' non-degenerate triangles.

\begin{Theorem} \label{irrational_theorem}
	Let $a \le b \le c$ be positive reals such that $c<a+b$ and $\frac{a}{c}$ or $\frac{b}{c}$ is irrational. Put $T=T(a,b,c)$.~If \hypertarget{star}{}
	\begin{itemize}[noitemsep,topsep=0pt]
		\item[$(\bigstar)$] $a = q_1 \xi$, $b = q_2 \eta$, $c = p_1 \xi + p_2 \eta$ for some odd integers $p_1, p_2, q_1, q_2$ and reals $\xi, \eta$ such that $\frac{\xi}{\eta}$ is irrational,
	\end{itemize}
	then $\chi(\mathbb{R}_\infty^2, T) = 3$. Otherwise, $\chi(\mathbb{R}_\infty^2, T) = 2$.
\end{Theorem}


For a degenerate triangle $T(a,b,a+b)$, several results, including the next one,  follow from~\cite{frankl2021max}, where much more general problems were studied. We provide its short proof here for completeness.

\begin{Proposition}
    Let $a, b$ be either coprime positive integers such that $a \equiv b \pmod 3$ or positive reals such that $\frac{a}{b}$ is irrational. Then $\chi \big(\mathbb{R}_\infty^2, T(a, b, a + b)\big) = 3$.
\end{Proposition}

\begin{proof}
	For the lower bound, observe that every five-point subset of a nine-element set $\{0,a,a+b\}^2 \subset \R^2$ contains an $\ell_\infty$-isometric copy of $T = T(a, b, a + b)$ and thus $\chi(\mathbb{R}_\infty^2, T) \ge \lceil\frac{9}{4}\rceil=3$ (cf. \cite[Proposition~17]{frankl2021max}).
	
	For the matching upper bound, it is crucial to note that at least one of the two projections of an $\ell_\infty$-isometric copy of $T$ onto basic axes is a translation of either $\{0,a,a+b\}$ or it's `reflection' $\{0,b,a+b\}$, as one would expect. Indeed, the largest distance $a+b$ between two of the vertices is determined by the absolute value of the difference of their either $x$- or $y$-coordinates. Then it is not hard to see that the projection of the third vertex onto the same axis is uniquely determined by lengths of the remaining two edges, i.e., which one is $a$, and which one is $b$ (cf. \cite[Lemma~12]{frankl2021max}).
	
	Now we claim that the coloring in the right-hand size of \Cref{fig1} does not contain a monochromatic $\ell_\infty$-isometric copy of $T$ if  $a$ and $b$ are coprime positive integers such that $a \equiv b \pmod 3$. Indeed, the projection of every color class onto each basic axis is a translation of $A = \{ x \in \mathbb{R} : \lfloor x \rfloor \not\equiv 0 \pmod 3 \}$. However, it is clear that $A$ does not contains translates of both $\{0,a,a+b\}$ and $\{0,b,a+b\}$ (cf. \cite[Corollary~20]{frankl2021max}).
	
	The case of irrational $\frac{a}{b}$ is similar. Let $H\subset \R$ be a \textit{Hamel basis} such that $a,b \in H$, see \Cref{sect_irr}. For $z \in \R$, $h \in H$, let $\lambda_h(z)$ be the coefficient of $h$ in the unique representation of $z$ as a finite linear combination of the elements of $H$ with rational coefficients. Now we color each point $\z=(x,y)\in \R^2$ according to the color of $\big(\lambda_a(x)+\lambda_b(x), \lambda_a(y)+\lambda_b(y)\big) \in \Q^2$ in the right-hand size of \Cref{fig1}. Arguing as in the previous paragraph, we can show that no $\ell_\infty$-isometric copy of $T$ is monochromatic (cf. \cite[Proposition~24]{frankl2021max}).
\end{proof}

In the remaining case $a \not\equiv b \pmod{3}$, we conjecture that three colors are not enough, which we verified for $a+b\le 7$ by computer search\footnote{We observed that every $3$-coloring of rectangular grids of sizes $4\times 5$, $6 \times 6$, $8\times 9$, $10\times 10$, $12\times 12$ and $12\times 13$ contains a monochromatic $\ell_\infty$-isometric copy of $T$ for $T$ = $T(1, 2, 3)$, $T(1, 3, 4)$, $T(2, 3, 5)$, $T(1, 5, 6)$, $T(1, 6, 7)$ and $T(3, 4, 7)$, respectively.}.

\begin{Conjecture} \label{conj}
	Let $a,b$ be coprime positive integers such that $a \not\equiv b \pmod{3}$. Then $\chi\big(\R_\infty^2,T(a,b,a+b)\big)=4$.
\end{Conjecture}

In what follows, for a shorthand, we refer to \textit{$\ell_\infty$-distances} and \textit{$\ell_\infty$-isometric copies} simply as \textit{distances} and \textit{copies}, respectively. Whenever we consider a two-coloring of the plane, we call these colors red and blue or $0$ and $1$ for clarity.

%% file: prelim.tex
In this section, we assume that the side lengths $a, b, c$ of a triangle $T = T(a, b, c)$ are reals satisfying $a \leq b \leq c$ and $c< a + b$, i.e., that $T$ is non-degenerate. First, we give a necessary condition for a triple of points to form a copy of $T$, and then we find some properties satisfied by every two-coloring of the plane containing no monochromatic copies of $T$.

\subsection{Copies of a non-degenerate triangle}

\begin{Lemma}\label{lemma_1}
	Let $\z_1=(x_1,y_1), \z_2=(x_2,y_2), \z_3=(x_3,y_3) \in \R^2$ be a copy of $T$. Then at least one of the differences $|y_1-y_2|, |y_2-y_3|, |y_3-y_1|$ equals either $a$, or $b$, or $c$. Moreover, at least one of the differences $|x_1+y_1-x_2-y_2|, |x_2+y_2-x_3-y_3|, |x_3+y_3-x_1-y_1|$ equals either $a+b-c$, or $c+a-b$, or $b+c-a$.
\end{Lemma}

\begin{proof}
	The distance between two points is determined by the absolute value of the difference of their either $x$- or $y$-coordinates. Therefore, one of the axes, say $x$, determines at least two of the distances $\|\z_1-\z_2\|_\infty, \|\z_2-\z_3\|_\infty, \|\z_3-\z_1\|_\infty$ by the pigeonhole principle. Assume that $\|\z_2-\z_3\|_\infty = |x_2-x_3|=a$, $\|\z_3-\z_1\|_\infty = |x_3-x_1|=b$. Observe that $x_3$ cannot lie between $x_1$ and $x_2$, since in that case we would get that $c = \|\z_1-\z_2\|_\infty \ge |x_1-x_2|=a+b$, a contradiction. Hence, let us assume that $x_1=x_3+b$, $x_2=x_3+a$. This implies that $|x_1-x_2|=b-a<c = \|\z_1-\z_2\|_\infty$ and thus $|y_1-y_2|=c$ as desired. To prove the second half of the statement, note that if $y_1=y_2+c$, then $|x_1+y_1-x_2-y_2|=b+c-a$, while if $y_1=y_2-c$, then $|x_1+y_1-x_2-y_2|=c+a-b$ as desired. The same reasoning works as well in all the remaining cases, corresponding to possible permutations of axes, indexes, and side lengths.  
\end{proof}

This simple statement immediately gives the following sufficient condition for a \textit{horizontal coloring}, which is constant on every horizontal line $y=y_0$, or for a \textit{diagonal coloring}, which is constant on every diagonal line $x+y=y_0$, to contain no monochromatic copies of $T$.

\begin{Corollary}\label{lemma_2}
	The following two statements are valid:
	\begin{enumerate}
		\item Let $\bar{C}(\cdot)$ be a coloring of the line such that no two points at distance $a$, $b$, or $c$ apart are monochromatic. Then the corresponding horizontal coloring of the plane, defined by the equation $C(x,y) = \bar{C}(y)$, contains no monochromatic copies of $T$.
		\item Let $C'(\cdot)$ be a coloring of the line such that no two points at distance $a+b-c$, $c+a-b$, or $b+c-a$ apart are monochromatic. Then the corresponding diagonal coloring of the plane, defined by the equation $C(x,y) = C'(x+y)$, contains no monochromatic copies of $T$.
	\end{enumerate}
\end{Corollary}

Perhaps surprisingly, later we will show that one can consider only these two types of colorings to find the exact value of $\chi(\R_\infty^2, T)$ for non-degenerate triangles. More formally, we will prove that for every coloring of the plane containing no monochromatic copies of $T$, there exists either a horizontal or a diagonal coloring that uses the same number of colors and also contains no monochromatic copies of $T$.

\subsection{Patterns in an arbitrary two-coloring}

For this subsection, let us fix a red-blue coloring of the plane such that no copy of $T$ is monochromatic. We call a vector $(x_0,y_0)$ its \textit{period} (resp. \textit{anti-period}) if for all $x,y\in \R$, the colors of two points $(x,y)$ and $(x+x_0,y+y_0)$ are the same (resp. distinct). It is clear that the addition of these vectors resembles the multiplication of signs: the sum of two periods or two anti-periods is a period, while the sum of an anti-period and a period is an anti-period.

\begin{Lemma} \label{anti-periods}
	The following three statements are valid:
	\begin{enumerate}
		\item If one of the four vectors $(\pm a, c-b)$, $(c-b, \pm a)$ is not an anti-period, then there exists a monochromatic axis-parallel segment\footnote{To clarify a possible ambiguity here and in what follows, by a \textit{segment} we mean a closed line segment.} of length $b+c-a$.
		\item If one of the four vectors $(\pm b, a-c)$, $(a-c, \pm b)$ is not an anti-period, then there exists a monochromatic axis-parallel segment of length $c+a-b$.
		\item If one of the four vectors $(\pm c, b-a)$, $(b-a, \pm c)$ is not an anti-period, then there exists a monochromatic axis-parallel segment of length $a+b-c$.
	\end{enumerate}
\end{Lemma} 
\begin{proof}
	If $(b-a, c)$ is not an anti-period, then there exist $x_1,y_1 \in \R$ such that the two points $(x_1,y_1)$ and $(x_1+a-b,y_1-c)$ are of the same color, say, both are red. It is easy to check that this pair together with an arbitrary point from the segment $\{(x_1-b, y): y_1-b \le y \le y_1+a-c\}$ form a copy of $T$. Since there are no red copies of $T$, we conclude that this vertical segment of length $a+b-c$ is entirely blue, as desired. Similar arguments work for the remaining $11$ vectors as well.
\end{proof}

\begin{Corollary} \label{anti-periods-2}
	If no axis-parallel segment of length $c+a-b$ is monochromatic, then all eight vectors $(\pm a, c-b)$, $(c-b, \pm a)$, $(\pm b, a-c)$, $(a-c, \pm b)$ are anti-periods, and all six\footnote{Some of these vectors may coincide if $a=b$, or $b-c$.} vectors $(2a,0)$, $(2b,0)$, $(2c,0)$, $(0,2a)$, $(0,2b)$, $(0,2c)$ are periods. Moreover, if there also exist $n,m,k \in \Z$ such that $0<2an+2bm+2ck\le a+b-c$, then four vectors $(\pm c, b-a)$, $(b-a, \pm c)$ are anti-periods as well.
\end{Corollary}
\begin{proof}
	Since $b+c-a \ge c+a-b$, the condition also implies that no axis-parallel segment of length $b+c-a$ is monochromatic. Now the first and the second parts of \Cref{anti-periods} give us eight desired anti-periods. To see that $(2a,0)$, $(2b,0)$, $(2c,0)$ are periods, we observe that
	$(2a,0) = (a,c-b)-(-a,c-b)$, $(2b,0) = (b,a-c)-(-b,a-c)$ $(2c,0)=(c-b,a)+(c-b,-a)+(2b,0)$. The case of $(0,2a)$, $(0,2b)$, $(0,2c)$ is similar.
	
	To prove the second half of the statement, note that both non-zero vectors $(2an+2bm+2ck, 0)$ and $(0, 2an+2bm+2ck)$ are periods as they are linear combinations of the other ones. Furthermore, observe that the inequality $2an+2bm+2ck \le a+b-c$ implies that no axis-parallel segment of length $a+b-c$ is monochromatic. Indeed, otherwise, the line containing this segment would also be monochromatic by the periodicity, which is a contradiction. Now the third part of \Cref{anti-periods} gives us four remaining anti-periods.
\end{proof}

\begin{Lemma}\label{lemma_3}
If $a<b$ and there exists a monochromatic axis-parallel segment of length $c+a-b$
, then there also exists a monochromatic axis-parallel line.
\end{Lemma}

\begin{proof}
	First, we consider a simple but illustrative special case when $b=c$. Let us assume without loss of generality that the monochromatic  axis-parallel segment of length $a=c+a-b$ is $I \coloneqq \{(x,0): 0\le x \le a\}$ and that it is red. It is clear that its two  endpoints $(0,0)$ and $(a,0)$ form a copy of $T$ together with an arbitrary point of the segment $J\coloneqq \{(x,b): a-b\le x \le b\}$ of length $2b-a=a+2(b-a)>a$. Since there are no red copies of $T$, we conclude that $J$ is entirely blue. Applying a similar argument to each pair of blue points from $J$ at distance $a$ apart, we get that the segment $I' \coloneqq \{(x,0): 2a-2b\le x \le 2b-a\}$ of length $a+4(b-a)$ is entirely red, see \Cref{fig3}. Proceeding in the same manner, we find a family of the concentric red segments $I\subset I' \subset I'' \subset \dots$ on the line $y=0$ whose lengths form an increasing arithmetic progression with the common difference of $2(b-a)$. Therefore, the horizontal line $y=0$ is entirely red.
	
	\begin{figure}[h]
		\includegraphics[scale=0.25]{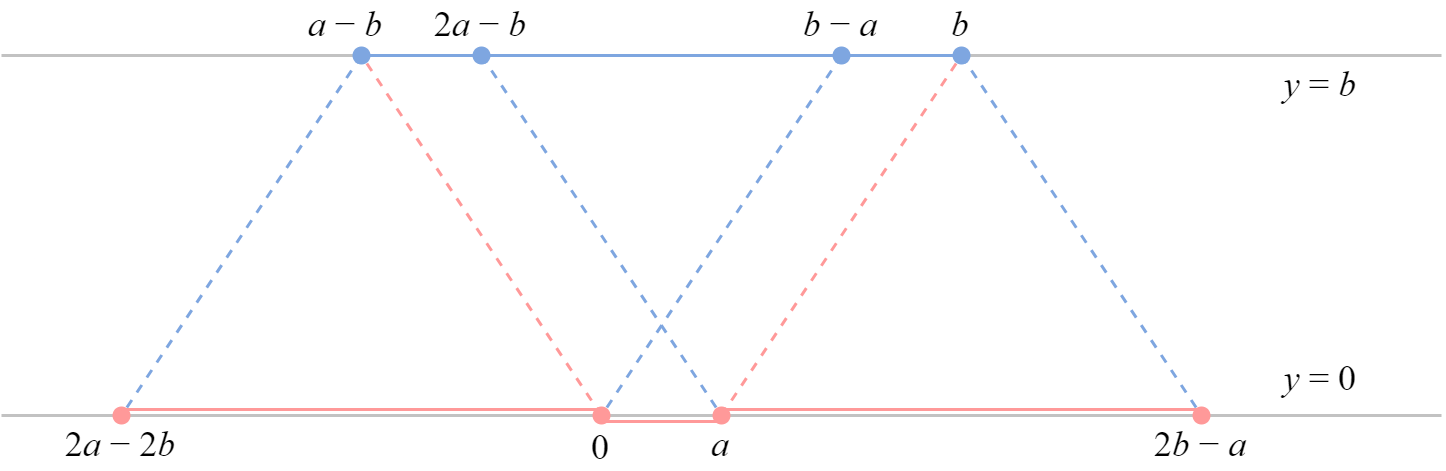}
		\centering
		\caption{Extension of a red segment in the `isosceles' case.}
		\label{fig3}
	\end{figure}
	
	Now we consider the remaining general case when $b<c$. If all four vectors $(\pm a, c-b)$, $(c-b, \pm a)$ are anti-periods, then both $(2a,0) = (a,c-b)-(-a,c-b)$ and $(0, 2c-2b) = (a,c-b)+(-a,c-b)$ are periods. Similarly, $(0,2a)$ and $(2c-2b,0)$ are also periods. It is clear that one of the two positive numbers $2a, 2c-2b$, which are the lengths of these periods, does not exceed their half-sum $c+a-b$, which is the length of the monochromatic axis-parallel segment given by the condition. Hence, the line containing this segment is also monochromatic due to the periodicity, and we are done.
	
	Otherwise, one of the four vectors $(\pm a, c-b)$, $(c-b, \pm a)$ is not an anti-period, and the first part of \Cref{anti-periods} ensures the existence of a longer monochromatic axis-parallel segment $I$ of length $b+c-a > c+a-b$. Let us assume without loss of generality that $I=\{(x,0): 0\le x \le b+c-a\}$ and that it is red. Suppose that we could show that these assumptions yield that the segment $J \coloneqq \{(x,a): b-c \le x \le 2c-a\}$ of length $b+c-a+2(c-b)>b+c-a$ is entirely blue. Then, by applying this argument repeatedly as earlier in the `isosceles' case, we would prove that each of the concentric segments $I\subset I' \subset I'' \subset \dots$ on the line $y=0$ whose lengths form an increasing arithmetic progression with the common difference of $2(c-b)$ is red, and thus the horizontal line $y=0$ is also entirely red, as desired.
	
	To show that $J$ is blue, we consider it as the union of four shorter segments first:
	\begin{align*}
		J'_1 \coloneqq&\ \{(x,a): b-c \le x \le b-a\}, \hspace{2mm} J'_2 \coloneqq \{(x,a): b-a \le x \le 2b-a\}, \\
		J_1 \coloneqq&\ \{(x,a): c \le x \le 2c-a\}, \hspace{6mm} J_2 \coloneqq \{(x,a): c-b \le x \le c\}.
	\end{align*}
	Indeed, to see that $J = J'_1\cup J'_2 \cup J_2 \cup J_1$, it is enough to check that $J'_2$ intersects $J_2$, i.e., that $2b-a\ge c-b$, which instantly follows from the triangle inequality $c\le a+b$. Moreover, observe that the images of $J'_1$ and $J'_2$ under the reflection about the vertical line passing through the midpoint of the segment $I$, namely about the line $2x=b+c-a$, are $J_1$ and $J_2$, respectively, see \Cref{fig4}. Thus, to complete the proof, it is sufficient to show that the latter two segments are blue.
	
	\begin{figure}[h]
		\includegraphics[scale=0.17]{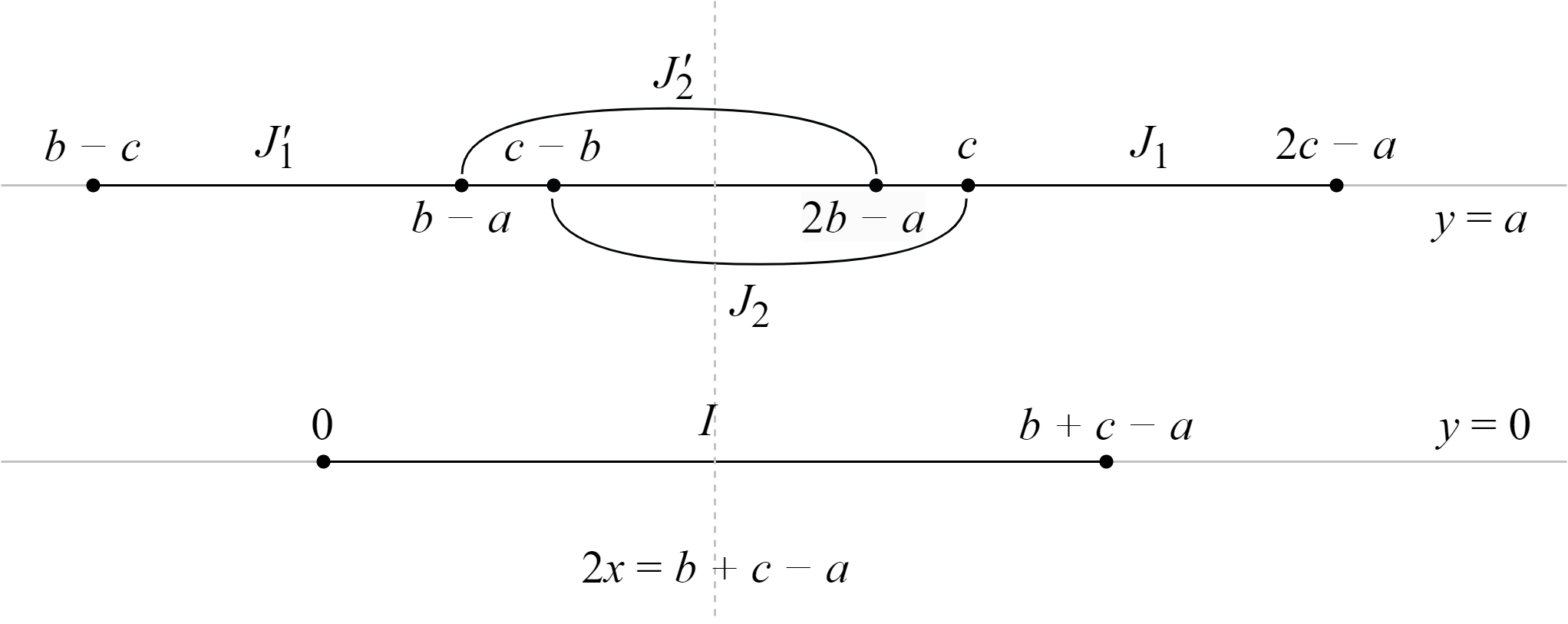}
		\centering
		\caption{A partition of $J$ into four shorter segments and its axis of symmetry.}
		\label{fig4}
	\end{figure}
	
	The case of $J_1$ is simple. For each $(x_0,a)\in J_1$, it is easy to check that both points $(x_0-c,0)$ and $(x_0+b-c,0)$ lie on $I$, and so they are red. Besides, this triple forms a copy of $T$, and thus $(x_0,a)$ is blue.

    The case of $J_2$ is trickier. Fix an arbitrary $\z_1=(x_1,a) \in J_2$ and suppose it is red. Consider the following four points: $\z_2=(x_1+b-c, 0), \z_3=(x_1+b, c+a-b), \z_4=(x_1+b,a-b), \z_5=(x_1+b-a,a)$, see \Cref{fig5}. Note that $\z_2 \in I$ and so it is red. Both triples $\z_1,\z_2,\z_3$ and $\z_1,\z_2,\z_4$ are copies of $T$ containing two red points $\z_1, \z_2$, and thus both $\z_3$ and $\z_4$ are blue. Finally, since $\z_3, \z_4, \z_5$ is also a copy of $T$, the point $\z_5$ is red. Observe that the distance between $\z_1$ and $\z_5$ is less than the length of $J_1$, namely that $b-a < c-a$. Therefore, we have $\z_5 \in J_2\cup J_1$. If the red point $\z_5$ belongs to $J_1$, then we immediately get a contradiction, since this segment is entirely blue as we have shown earlier. Otherwise, we apply a similar argument with $\z_5$ playing the role of $\z_1'$ to show that the point $\z_5' = (x_1+2b-2a,a)$ is also red. Proceeding in the same manner, we obtain a sequence of the red points $\z_1, \z_1', \z_1'', \dots$ on the horizontal line $y=a$, whose $x$-coordinates form an increasing arithmetic progression with the common difference of $b-a$. Once a point of this sequence falls onto the blue segment $J_1$, we terminate the process and get the desired contradiction.   
\end{proof}

\begin{figure}[h]
	\includegraphics[scale=0.19]{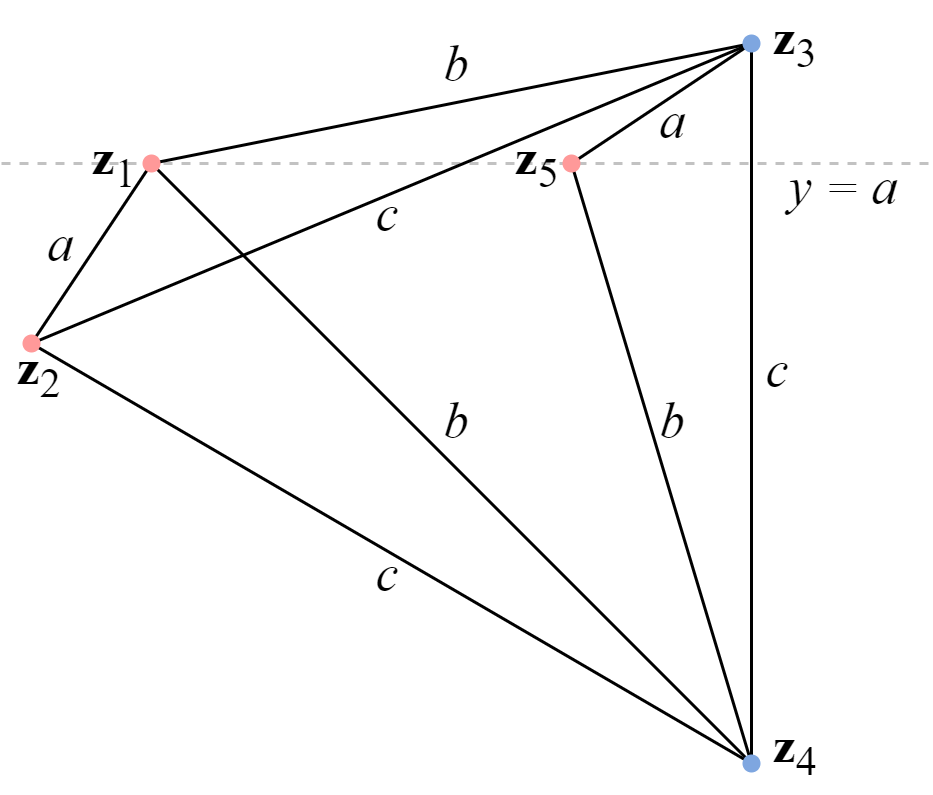}
	\centering
	\caption{Five points determine three copies of $T$.}
	\label{fig5}
\end{figure}

\begin{Lemma} \label{lines}
	If the horizontal line $y=0$ is red, then both lines $y=a$ and $y=b$ are blue. Moreover, if there also exist $n,m \in \Z$ such that $n+m$ is even and $c-b \le an+bm \le a$, then the line $y=c$ is blue as well.
\end{Lemma}
\begin{proof}
	Each point $(x_0,a)$ on the line $y=a$ forms a copy of $T$ together with two red points $(x_0-c,0)$ and $(x_0+b-c,0)$. Thus the line $y=a$ is entirely blue. Similarly, each point $(x_0,b)$ on the line $y=b$ forms a copy of $T$ together with two red points $(x_0-c,0)$ and $(x_0+a-c,0)$. Thus the line $y=b$ is also entirely blue.
	
	To prove the second half of the statement, observe that the first half implies that the color of the horizontal line $y=an+bm$ is determined by the parity of $n+m$. In particular, if $n+m$ is even, then this line is red. Now it is easy to see that each point $(x_0,c)$ on the line $y=c$ forms a copy of $T$ together with two red points $(x_0-b,an+bm)$ and $(x_0+a-b,0)$. Hence, the line $y=c$ is also entirely blue, as desired.
\end{proof}

%% file: proof.tex
First of all, we note that for every non-degenerate triangle $T=T(a,b,c)$ with integer side lengths $a\le b \le c$ such that $\gcd(a,b,c)=1$, the upper bound $\chi(\R_\infty^2,T) \le 3$ is immediate from the first half of \Cref{lemma_2} and the following special case\footnote{One can also deduce this result from~\cite[Theorem~3]{heuberger2003planarity}. Indeed, it is sufficient to color $\Z_n$ for some $n$ and then extend it to the whole line by coloring each $x\in\R$ with the color of $\lfloor x \rfloor \bmod n$. By taking $n$ equal to $a+b$, $b+c$, or $c+a$, we can identify two of the three forbidden distances, namely $a$ and $b$, $b$ and $c$, or $c$ and $a$, respectively. A direct case analysis shows that in at least one of these three cases, there exists a desired coloring of $\Z_n$ in only three colors.} of \cite[Corrolary~2.1]{Zhu} due to Zhu.

\begin{Theorem} \label{ZhuTh}
	Let $a\le b\le c$ be positive integers such that $c<a+b$ and $\gcd(a, b, c)=1$. Then there exists a three-coloring of the line such that no two points at distance $a, b,$ or $c$ apart are monochromatic.
\end{Theorem}

Therefore, to complete the proof, we only need to show that there exists a two-coloring of the plane such that no copy of $T$ is monochromatic if and only if either \refheart\ or \refdiamond\ holds. Let us begin by showing the sufficiency of these conditions using the following two explicit colorings.

\begin{Proposition} \label{casei}
	In the notation of \Cref{main_theorem}, assume that \refheart\ holds. Consider a diagonal two-coloring of the plane defined by $C(x,y) = \lfloor x+y\rfloor \bmod 2$. Then no copy of $T$ is monochromatic.
\end{Proposition}
\begin{proof}
	It is clear that $\lfloor x_1 \rfloor \not\equiv \lfloor x_2 \rfloor \pmod{2}$ whenever $x_1,x_2 \in \R$ are at odd distance apart. Since all three values $a+b-c$, $c+a-b$, and $b+c-a$ are odd by \refheart, the second half of \Cref{lemma_2} completes the proof.
\end{proof}

\begin{Proposition} \label{caseii}
	In the notation of \Cref{main_theorem}, assume that \refdiamond\ holds. Consider a horizontal two-coloring of the plane defined by $C(x,y) = \lfloor y/d\rfloor \bmod 2$, where $d=\gcd(a,b)$. Then no copy of $T$ is monochromatic. 
\end{Proposition}
\begin{proof}
	Assume the contrary, namely that for some $\z_1=(x_1, y_1), \z_2=(x_2, y_2), \z_3=(x_3, y_3) \in \R^2$ that form a copy of $T$, the values $\lfloor y_1/d \rfloor, \lfloor y_2/d \rfloor, \lfloor y_3/d \rfloor$ are of the same parity, say all three are even. By \Cref{lemma_1}, we can assume without loss of generality that $y_1 - y_2$ equals either $a$, or  $b$,  or $c$. Note that the former two cases immediately yield a contradiction since both fractions $a/d$ and $b/d$ are odd by \refdiamond. So in what follows we suppose that $y_1-y_2=c$. In particular, this implies that $\|\z_1-\z_2\|_\infty = c$, and thus one of the distances $\|\z_2-\z_3\|_\infty$, $\|\z_3-\z_1\|_\infty$ equals $a$, while the other one equals $b$. 
	
	It is easy to check that if $\|\z_2-\z_3\|_\infty =a$, $\|\z_3-\z_1\|_\infty =b$, then $y_1-b \le y_3 \le y_2+a$. Observe that both $\lfloor(y_1-b)/d\rfloor = \lfloor y_1/d\rfloor - b/d$ and $\lfloor(y_2+a)/d\rfloor = \lfloor y_2/d\rfloor + a/d$ are odd. Moreover, the difference $(y_2+a)-(y_1-b)=a+b-c$ does not exceed $d$ by \refdiamond. Therefore, the parity of $\lfloor y/b \rfloor$ can not change from odd to even and back again as $y$ ranges between these two extreme values. Hence, $\lfloor y_3/d \rfloor$ is also odd, and we see the contradiction. In the remaining case when $\|\z_2-\z_3\|_\infty =b$, $\|\z_3-\z_1\|_\infty =a$, we have $y_1-a \le y_3 \le y_2+b$, and a similar argument completes the proof.
\end{proof}

\noindent
To prove the second half of the theorem, let us observe that if neither \refheart\ nor \refdiamond\ holds, then either\hypertarget{club}{}
\begin{itemize}[noitemsep,topsep=0pt]
	\item[$(\clubsuit)$] $a$ and $b$ are of different parity, $c$ is odd, or \hypertarget{spade}{}
	\item[$(\spadesuit)$]  $a$ and $b$ are odd, $c$ is even, and $c< a+b-\gcd(a,b)$.
\end{itemize}
So it remains only to show that in each of these two cases, there are no two-colorings of the plane with no monochromatic copies of $T$. Let us assume the contrary and fix an arbitrary such red-blue coloring.

First, we suppose that no axis-parallel segment of length $c+a-b$ is monochromatic. On the one hand, the first half of \Cref{anti-periods-2} implies that all six vectors $(2a,0)$, $(2b,0)$, $(2c,0)$, $(0,2a)$, $(0,2b)$, $(0,2c)$ are periods, and so are all their linear combinations. Since $\gcd(2a,2b,2c)=2$, we conclude that both $(2,0)$ and $(0,2)$ are periods. Hence, every vector such that both its coordinates are even integers is also a period. On the other hand, note that $a+b-c$ is a positive even integer, and thus $a+b-c\ge 2 = \gcd(2a,2b,2c)$. Therefore, we can also apply the second half of \Cref{anti-periods-2} in our case to find twelve anti-periods in total including $(a, c-b)$, $(b, a-c)$ and $(c, b-a)$. However, both coordinates of one of these three vectors are even integers, and so this vector should be a period instead, a contradiction.

Second, we suppose that there exists a monochromatic axis-parallel segment of length $c+a-b$. Besides, note that each of the conditions \refclub\ and \refspade\ yields that $a<b$. 
So we can apply \Cref{lemma_3} to find a monochromatic axis-parallel line. Without loss of generality, let us assume that this line is given by the equation $y=0$ and that it is entirely red. Now it is easy to deduce from the first half of \Cref{lines} that for all $i,j \in \Z$ such that $i+j$ is odd, the horizontal line $y=ai+bj$ is blue. If \refclub\ holds, then we obtain a contradiction by taking $i=b$, $j=-a$.

If \refspade\ holds, we use a slightly more complex argument. Observe that there exist $n,m \in \Z$ such that $an+bm=a-\gcd(a,b)$. Since both $a$ and $b$ are odd, we conclude that $n+m$ is even. Moreover, the inequality $c-b\le an+bm = a-\gcd(a,b) < a$ is immediate from \refspade. Therefore, we can also apply the second half of \Cref{lines} in our case to deduce that the horizontal line $y=ai+cj$ is blue for all  $i,j \in \Z$ such that $i+j$ is odd. Finally, we obtain the desired contradiction by taking $i=c$, $j=-a$.

%% file: proof2.tex
Our arguments in this section heavily rely on the axiom of choice and, in particular, on the notion of a \textit{Hamel basis}, see, e.g., \cite[Chapter~4]{KuczmaMarek}. Recall that a subset $H \subset \R$ is a \textit{Hamel basis} if every $x\in \R$ may be uniquely represented as a finite linear combination of the elements of $H$ with rational coefficients. We denote the coefficient of $h\in H$ in this representation by $\lambda_h(x)$. Note that for every $h \in H$, we have $\lambda_h(x+y)=\lambda_h(x)+\lambda_h(y)$ for all $x,y\in \R$. Let us also recall that every set of linearly independent (over $\Q$) reals can be completed to a Hamel basis, see \cite[Theorem~4.2.1]{KuczmaMarek}.

Let $T=T(a,b,c)$ be a non-degenerate triangle with side lengths $a\le b \le c$ such that at least one of the fractions $a/c$, $b/c$ is irrational. As in the previous section, we begin with an upper bound on $\chi(\R_\infty^2,T)$, which is immediate from the second half of \Cref{lemma_2} and the following statement. 

\begin{Proposition} \label{irr_upper}
	In the notation of \Cref{irrational_theorem}, there exists a three-coloring of the line such that no two points at distance from $D\coloneqq\{a+b-c, c+a-b, b+c-a\}$ apart are monochromatic. Moreover, if \refstar\ does not hold, then there even exists a two-coloring of the line with the same property. 
\end{Proposition}
\begin{proof}
	First, we assume that $c$ cannot be written as a linear combination of $a$ and $b$ with rational coefficients.
	
	We claim that there exists a Hamel basis $H \subset \R$ such that $c \in H$ and $\lambda_c(a) =\lambda_c(b)=0$. To show this, let us consider two cases. If the set $\{a,b,c\}$ is linearly independent over $\Q$, then every Hamel basis $H$ such that $\{a,b,c\}\subset H$ suffices. Indeed, $a=a$ and $b=b$ are the unique representation of $a$ and $b$, respectively, and thus $\lambda_c(a)=\lambda_c(b)=0$. Otherwise, there exist $r_1,r_2,r_3 \in \Q$ such that $r_1a+r_2b+r_3c=0$. It is clear that $r_3=0$ by our assumption. Then every Hamel basis $H$ such that $\{a,c\}\subset H$ suffices, since $a=a$ and $b=-\frac{r_1}{r_2}a$ are the unique representations of $a$ and $b$, respectively, and thus  $\lambda_c(a)=\lambda_c(b)=0$, as claimed. \vspace*{0.5mm}
	
	With this Hamel basis, it is easy to see that for all $x_1,x_2 \in \R$ such that $|x_1-x_2| \in D$, we have $|\lambda_c(x_1)-\lambda_c(x_2)|= 1$. In particular, this implies that $\lfloor \lambda_c(x_1) \rfloor \not\equiv \lfloor \lambda_c(x_2) \rfloor \pmod{2}$, and thus the equation $C'(x) = \lfloor \lambda_c(x)\rfloor \bmod 2$ defines the desired two-coloring of the line.
	
	Now we proceed to the other alternative in which $c$ can be written as a linear combination of $a$ and $b$ with rational coefficients, namely $c=\frac{p_1}{q_1}a+\frac{p_2}{q_2}b$ for some integers $p_1,q_1,p_2,q_2$ such that $\gcd(p_1, q_1) = \gcd(p_2, q_2) = 1$. Taking $\xi\coloneqq a/q_1, \eta \coloneqq b/q_2$, we rewrite $c=p_1\xi+p_2\eta$. Moreover, note that $\xi/\eta$ is irrational, since otherwise both fractions $a/c$ and $b/c$ are rationals. Let $H \subset \R$ be a Hamel basis such that $\xi,\eta \in H$. It is easy to check that for all $x_1,x_2 \in \R$ such that $|x_1-x_2| \in D$, we have $|\lambda_\xi(x_1)-\lambda_\xi(x_2)| \in D_1 \coloneqq \{|p_1+q_1|, |p_1-q_1|\}$, and $|\lambda_\eta(x_1)-\lambda_\eta(x_2)| \in D_2 \coloneqq \{|p_2+q_2|, |p_2-q_2|\}$. This observation implies that if $C_1'(\cdot)$ colors the integers such that no two of them at distance from $D_1$ apart are monochromatic, then $C_1'(\lfloor \lambda_\xi(x_1)\rfloor) \neq C_1'(\lfloor \lambda_\xi(x_2)\rfloor)$ for all $x_1,x_2 \in \R$ such that $|x_1-x_2| \in D$. Similarly, if $C_2'(\cdot)$ colors the integers such that no two of them at distance from $D_2$ apart are monochromatic, then $C_2'(\lfloor \lambda_\eta(x_1)\rfloor) \neq C_1'(\lfloor \lambda_\eta(x_2)\rfloor)$ for all $x_1,x_2 \in \R$ such that $|x_1-x_2| \in D$. So in the rest of this proof, we consider only these `integral' colorings.
	
	Note that if one of the two coprime integers $p_1,q_1$ is even, then the other one is odd. The same is true for the other pair $p_2,q_2$ as well. Therefore, if \refstar\ does not hold, namely if at least one of the four integers $p_1,q_1,p_2,q_2$ is even, then $D_i$ contains only odd integers for some $i=1,2$. Thus the two-coloring defined by $C_i'(n) = n \bmod{2}$ contains no monochromatic integers at distance form $D_i$ apart, as desired.
	
	Finally, it remains only to show that if \refstar\ holds, then for some $i=1,2$, there exists a three-coloring of the integers such that no two of them at distance from $D_i$ apart are monochromatic. Roughly speaking, this is true because the chromatic number of each \textit{$2$-degenerate} graph does not exceed $3$. More formally, we order the integers such that their absolute values do not decrease and build the desired coloring iteratively. Assume that after $m$ steps, we have already colored the first $m$ integers, and consider the next one, say $n$. Among the four integers of the form $n\pm d$, $d \in D_i$, at least two have absolute values larger than $|n|$, and thus they have not been colored yet. Hence, there are at most two colors that we cannot use to color $n$, and thus at least one color is available. Proceeding in the same manner, we obtain the desired three-coloring.
	
	Observe that the previous argument requires that $0 \notin D_i$. However, this is always the case either for $D_1$ or for $D_2$. Indeed, if $0 \in D_1\cap D_2$, then $p_1=\pm q_1$ and $p_2=\pm q_2$. Therefore, one of the four expressions of the form $\pm a \pm b$ is equal to $c$, which contradicts the chain of inequalities $0<a\le b \le c < a+b$. 
\end{proof}

Once we obtained the aforementioned upper bound on $\chi(\R_\infty^2,T)$, it remains only to get the matching lower bound to complete the proof of \Cref{irrational_theorem}. Namely, we need to show that $\chi(\R_\infty^2,T)>2$ if \refstar\ holds, i.e., if $a = q_1 \xi$, $b = q_2 \eta$, $c = p_1 \xi + p_2 \eta$ for some odd integers $p_1, p_2, q_1, q_2$ and reals $\xi, \eta$ such that $\xi/\eta$ is irrational. As in the previous section, we do it by contradiction. So let us assume the contrary and fix an arbitrary red-blue coloring of the plane with no monochromatic copies of $T$.

First, we suppose that there exists a monochromatic axis-parallel segment of length $c+a-b$. Since we clearly have $a\neq b$ in our case, \Cref{lemma_3} implies the existence of a monochromatic axis-parallel line. Without loss of generality, let us assume that this line is given by the equation $y=0$ and that it is entirely red. Observe that the set $\{an+bm: n,m \in \Z\}$ is dense since the fraction $a/b$ is irrational. In particular, this implies that we can use both parts of \Cref{lines} to deduce that the horizontal line $y=ai+bj+ck$ is blue for all  $i,j,k \in \Z$ such that $i+j+k$ is odd. By taking three odd integers $i=p_1q_2$, $j=q_1p_2$, $k=-q_1q_2$, we get that the horizontal line $y=0$ is blue, a contradiction.

Second, we suppose that no axis-parallel segment of length $c+a-b$ is monochromatic. Recall that the set $\{an+bm: n,m \in \Z\}$ is dense, and thus we can use both parts of \Cref{anti-periods-2} to find twelve anti-periods of our coloring in total. In particular, their linear combination $\V_1\coloneqq (b+c-a,0) = (b-a,c)-(-a,c-b)-(a-c,b)$ is also an anti-period. Similar argument yields that both $\V_2\coloneqq (c+a-b,0)$ and $\V_3\coloneqq (a+b-c,0)$ are anti-periods as well. Consider three integers $i=(q_1p_2-q_1q_2)/2$, $j=(p_1q_2-q_1q_2)/2$ and $k=(p_1q_2+q_1p_2)/2$. On the one hand, their sum $i+j+k=p_1q_2+q_1p_2-q_1q_2$ is odd, and thus the linear combination $i\V_1+j\V_2+k\V_3$ is an anti-period. On the other hand, it is not hard to check that the latter vector equals $(0,0)$ which is a trivial period, a contradiction.

%% file: conclusions.tex
Though \Cref{conj} is the main open problem of our paper, let us briefly discuss some more questions.

\vspace{1mm}
\noindent \textbf{Finite domains with the same chromatic numbers.} Given a triangle $T$, the hypergraph version of the De Bruijn--Erdős theorem \cite{bruijn1951colour} implies that there exists a finite $V \subset \R^2$ with the following property. Whenever $V$ is colored in less then $\chi(\R_\infty^2,T)$ colors, it contains a monochromatic $\ell_\infty$-isometric copy of $T$. It is natural to ask what is the minimum size of such $V$. Some upper bound on this quantity in the non-degenerate case can be derived from our proofs here, but it would probably be far from optimal.

\vspace{1mm}
\noindent \textbf{Fractional chromatic numbers.} For a triangle $T$ with integer side lengths, it would be interesting to find the exact value of $d_f(\Z^2,T)$ defined as the maximum asymptotic density of a subset of the integer grid $\Z^2$ with no $\ell_\infty$-isometric copies of $T$. This is equivalent to determining the \textit{fractional chromatic number} $\chi_f(\R_\infty^2,T)$, see~\cite{chang1999distance}. In contrast to our previous results, this problem is open only for non-degenerate triangles. Indeed, for a degenerate triangle $T$, \cite[Proposition~17]{frankl2021max} yields the equality $d_f(\Z^2,T) = d_f(\Z,T)^2$ reducing this to a one-dimensional problem, which was later solved in~\cite{frankl2022solution}.

\vspace{1mm}
\noindent \textbf{Larger dimensions.} This situation is similar for multidimensional versions of these chromatic numbers. For every degenerate triangle $T$, the constant $c=c(T)>1$ such that $\chi(\R_\infty^n,T) = \big(c+o(1)\big)^n$ as $n\to \infty$ was found in~\cite{frankl2021max}. At the same time, such $c(T)$, if it exists, is known only for a few non-degenerate triangles $T$. For instance, assume that $T=T(a,a,b)$ is an isosceles triangle. If its base is not longer than its legs, i.e., if $b\le a$, then \cite[Theorem~5]{frankl2021max} implies that $\chi(\R_\infty^n,T) = \big(2+o(1)\big)^n$ as $n\to \infty$. However, if the base is longer, then we only know that $\big(3/2+o(1)\big)^n \le \chi(\R_\infty^n,T) \le \big(1+a/b+o(1)\big)^n$, where the upper bound is only slightly better than the trivial inequality $\chi(\R_\infty^n,T) \le \chi(\R_\infty^n) = 2^n$ when the triangle is close to be equilateral.

\vspace{1mm}
\noindent \textbf{Larger configurations.} Though some degenerate triangles require four colors to color the plane such that none of their $\ell_\infty$-isometric copies are monochromatic, for all sufficiently large forbidden configurations $\mathcal{M}$ three colors are always enough, see the proof of \cite[Proposition~7]{frankl2021max}. Perhaps, the condition $|\mathcal{M}|\ge 4$ is already sufficient for this. So the chromatic number $\chi(\R_\infty^2,\mathcal{M})$ is equal to either $2$ or $3$ for every large $\mathcal{M}$, and it is not hard to see that both options indeed take place, see \cite[Theorem~8]{frankl2021max}. Therefore, it is reasonable to ask if there is a simple criterion distinguishing them.

\vspace{1mm}
\noindent \textbf{Other metrics.} Counterparts of the Nelson's classic question on the chromatic number of the plane were studied for many other non-Euclidean metrics, see \cite{chilakamarri1991unit, exoo2021chromatic, geher2023note}. It might be interesting to study the colorings with no monochromatic copies of a given triangle for these metrics as well.

\subsection*{Acknowledgments}

\noindent
We thank N\'ora Frankl for the helpful discussion of \Cref{conj} and Panna Geh\'er for bringing the reference~\cite{aichholzer2019triangles} to our attention. 
The second author is supported by ERC Advanced Grant `GeoScape' No. 882971.